\documentclass[a4paper,review]{elsarticle}

\usepackage{lineno,hyperref}
\usepackage{pgfplots}
\pgfplotsset{compat=newest}
\nolinenumbers
\modulolinenumbers[0]

\journal{~}

\usepackage{tikz}
\usepackage{bm}
\usepackage{enumitem}
\usepackage{algorithm}
\usepackage{algorithmicx}
\usepackage{algpseudocode}
\usepackage{mathtools}
\usepackage{amssymb}
\allowdisplaybreaks

\begin{document}
	
	\begin{frontmatter}
		
		%\title{}
		%\tnotetext[mytitlenote]{}
		
		\title{Profit maximization via capacity control for distribution logistics problems}
		
		%\title{The capacity control problem for revenue maximization of distribution logistics systems}
		
		\author{Giovanni Giallombardo}
		\address{Dipartimento di Ingegneria Informatica, Modellistica, Elettronica e Sistemistica \\ Universit\`a della Calabria, 87036 Rende (CS), Italy \\Email: giallo@dimes.unical.it}

		\author{Francesca Guerriero}
		\address{Dipartimento di Ingegneria Meccanica, Energetica e Gestionale \\ Universit\`a della Calabria, 87036 Rende (CS), Italy \\ Email: francesca.guerriero@unical.it}
		
		\author{Giovanna Miglionico\corref{mycorrespondingauthor}}
		\address{Dipartimento di Ingegneria Informatica, Modellistica, Elettronica e Sistemistica \\ Universit\`a della Calabria, 87036 Rende (CS), Italy\\ Email: gmiglionico@dimes.unical.it}
%		\cortext[mycorrespondingauthor]{Corresponding author}
		%\ead{gmiglionico@dimes.unical.it}
		
		\begin{abstract}
			We consider a distribution logistics scenario where a shipping operator, managing a limited amount of resources, receives a stream of collection requests, issued by a set of customers along a booking time-horizon, that are referred to a future operational period. The shipping operator must then decide about accepting or rejecting each incoming request at the time it is issued, accounting for revenues, but also considering resource consumptions. In this context, the decision process is based on dynamically finding the best trade-off between the immediate return of accepting the request and the convenience of preserving capacity to possibly exploit more valuable future requests. We give a dynamic formulation of the problem aimed at maximizing the operator revenues, accounting also for the operational distribution costs.
			Due to the ``curse of dimensionality'', the dynamic program cannot be solved optimally. For this reason, we propose a mixed-integer linear programming approximation, whose exact or approximate solutions provide the relevant information to apply some commonplace revenue management policies in the real-time decision-making. Adopting a capacitated vehicle routing problem as an underlying distribution application, we analyze the computational behaviour of the proposed techniques on a set of academic test problems.

		\end{abstract}
		
		
	\end{frontmatter}
	
	%\linenumbers
	
	\section{Introduction}
	\label{section1}
	
	Revenue management (RM) techniques, originally developed in the airline industry in 1970s \cite{Belobaba87b}, with the aim of efficiently allocating fixed resource capacity to different classes of customers, have rapidly spread throughout other sectors of application, such as transportation \cite{guerriero12b,guerriero12a}, hotels \cite{goldman02,kimes11}, and restaurants \cite{guerriero14}.
	
	The growing relevance of RM is related to its practical importance in maximizing company revenue, while selling the same amount of products or services. In this context, customer segmentation, demand forecasting, pricing, and customer behavior play a crucial role \cite{talluri05}, and RM techniques can be very effective in dynamically estimating and controlling the unknown demand by appropriately setting prices, or by efficiently allocating resources.
	\color{black}In fact, along a discrete time horizon, \color{black}pricing techniques aim at dynamically setting different prices for the same product, while capacity-control techniques aim at dynamically allocating resource capacity on differentiated products/services that are generally booked in advance.
	Focusing, in particular, on the latter, \color{black}at each booking period \color{black}the main decision is to accept or reject an \color{black}incoming booking-request \color{black}for a given product/service. \color{black}Thus, such decision \color{black}is based on finding a trade-off between the immediate return of accepting the request, and the convenience of preserving capacity (i.e., resource availability) to possibly exploit more valuable \color{black}future \color{black}requests.

	RM techniques turn to be very effective when applied to industries with specific
	characteristics: a fixed capacity, a cost structure characterized by prevailing fixed costs, and a
	variable and segmented demand. The extension of such techniques to cases where some of such features are missing is not trivial. In fact, our aim is to deal with problems where, at the end of the \color{black}booking horizon, \color{black}non-negligible operational costs, which depend on the set of accepted booking-requests, arise. More precisely, it is assumed that the product/service provider has several operational options to satisfy all the accepted requests.

	%It is worth observing that in some applications, it is not possible to
	%neglect variable costs, since an accepted request produces both a revenue, obtained from selling the product/service,
	%and a cost variation due to the increased consumption of all the resources involved.
	%Moreover, the calculation of such cost variation typically requires a reconfiguration of the
	%whole  system, by means of the solution of a cost minimization operational problem.
	%Rather than dealing with revenue maximization, in such cases it is more appropriate to deal with a profit maximization
	%problem. This is the case of the problem addressed in this paper.
	The main focus of the article is to study RM capacity-control problems in the distribution logistics sector. In particular, we consider a logistic operator that collects parcels from a certain number of customers on a distribution network. At each time-period of the booking horizon when an additional service-request is received, in order to make an accept/reject decision, the operator must optimize a \color{black}profit \color{black}function that accounts for (i) the revenue associated to the current request, (ii) the \color{black}expected revenue associated to future requests, \color{black}and (iii) the final operational costs associated to servicing all the accepted requests. \color{black} As for the latter, we will restrict our attention to an operational problem of the vehicle-routing type, although the proposed approach is fairly general and can be suited, \textit{mutatis mutandis}, to different distribution problems. \color{black}
	%In fact, under such scheme, an additional request returning a negative
	%current profit might even be accepted in case the future profit is enough to balance the current
	%loss, this being the case when the additional request requires the opening of a new route where the
	%potential future customers are many.
		
	\color{black}We provide a dynamic formulation for the problem of accepting/rejecting the stream of incoming product-requests with the aim of maximizing the operator profit. \color{black} Due to ``the curse of dimensionality'', the dynamic programming
	model cannot be solved optimally. For this reason, in order to provide the logistic operator
	with a tool useful in making decisions, we propose an integer linear programming (ILP) approximation of the dynamic program, whose exact or approximate solutions provide the relevant information to apply some commonplace RM policies in the real-time decision-making.
	
	The problem is closely related to a class of RM techniques developed in the attended home delivery management. Such approaches are focussed on pricing techniques, that also account for customer choice behaviour, and typically rely on heuristics aimed at approximating opportunity costs \cite{Yang16, Yang17}.    
			
	\color{black} Our contribution is to apply RM capacity-control techniques to distribution logistics problems where operational costs can be optimized,  gaining a general insight into the structure of ILP approximations associated to the operational problem. In particular, for the underlying capacitated vehicle-routing problem (CVRP) considered here, we come up with a formulation that defines a novel subclass of CVRP whose aim is the profit maximization. \color{black}
	
	The rest of the paper is organized as follows. \color{black}In Sections \ref{section2.0} and \ref{section2}, respectively, we briefly survey some relevant literature regarding capacity-control problems, and review the standard methodological framework of multiple-resource capacity control (see \cite{bertsimas03,deboer02a,talluri05}). \color{black}In Section \ref{section3} we introduce a dynamic programming formulation for the capacity control in distribution logistics, and provide a specific development based on adopting a CVRP as an underlying distribution problem. In Section \ref{section4} we propose a deterministic approximation of the dynamic program, together with the related control policies. Numerical experiments  are presented in Section \ref{section5} and some conclusions are drawn in Section \ref{section6}.
	
	\section{Literature review}
	\label{section2.0}
	The fundamental capacity control decision is to accept or reject a booking request for a given product, based on finding a trade-off between the immediate return of accepting the request, and the future convenience of preserving capacity to exploit possibly more valuable requests. The underlying mathematical model is a dynamic and stochastic program. Among the first works in capacity control we recall \cite{Belobaba87b} where the seat inventory control in airlines is discussed, and \cite{Belobaba89a} where the expected marginal seat revenue model, taking into account the demand uncertainty, is presented and tested on real-life airlines data. Since the dynamic problem (DP) is impossible to solve in practice, several approximations of the capacity control problem have been developed and tested. In particular, we remember the approximations by deterministic linear
	programming (DLP) considered in \cite{deboer02a}, the  probabilistic nonlinear programming discussed in \cite{talluri05}, the randomized linear programming introduced in \cite{talluri99} and opportunity cost-based methods proposed in \cite{bertsimas03}. Relationships among  DP and DLP, and also the so-called asymptotic optimality properties, have been widely studied, see for some examples \cite{cooper02} and  \cite{talluri98}.
	
	Revenue management policies are  particularly efficient when applied to industries with specific characteristics: a fixed capacity, a cost structure characterized by large fixed costs and a variable and segmented demand. In the last few decades, the application of revenue management techniques has also been extended to many industries, we cite for example restaurants, cargo, hospitals, cruise lines, and golf. \color{black}Interesting reviews on revenue management features, techniques and sectors of application can be found in \cite{Chiang07,klein20}. \color{black}Recent papers applying revenue management policies to less conventional industries are \cite{guerriero12b,guerriero12a,guerriero14}.
	
	In order to widen the field of application of revenue management, some attempts have been done to drop typical assumptions. \color{black}The introduction of incompletely specified products \cite{gonsch20}, like flexible products \cite{gallego04} and opaque products \cite{Chen10}\color{black}has offered  the possibility of dealing with a less segmented and hence potentially larger demand, \color{black}as for such products industries can partially postpone resource-product assignment at the end of the booking horizon.  \color{black}
	Important attempts have also been made to relax some strong assumptions related to the customer behavior. The introduction of the choice model makes the purchase behavior dependent on current availability of products. A general analysis of the choice model based on dynamic programming is provided in \cite{talluri04a} for a single-leg setting. The extension to the \color{black}multiple-resource \color{black} setting is presented in \cite{vanRyzin04a,vanRyzin04b} and in \cite{MirandaBront09a}. In \cite{jiang14} another characterization of the customer choice behavior is considered, where a denied customer has the possibility of buying up to a product with different characteristics (prices, times or routes). 
	
	\color{black}Problems where operational decisions, depending on the accepted booking requests,  affect the revenue levels have already raised the interest of the revenue management community. In this respect the approaches proposed for managing the overbooking mechanism in the airline industry are among the most relevant examples.  \color{black}It is well known that  overbooking is implemented to handle the uncertainty in capacity utilization. For an example, tickets in the airline industry are often sold beyond the physical capacity of each aircraft with the aim of profitably exploiting cancellations and no-shows \cite{bertsimas03,Rothstein71,Subramanian99}. Depending on the latter ones, a penalty cost for serving extra requests is possibly incurred at the end of the booking horizon, given the potential need to bump some of the accepted requests.
	Other approaches, which account for such decision-dependent costs, deal with the resource-capacity reconfiguration, driven by revenue maximization, to handle demand variation. The reader is referred to \cite{Barnhart2009,Berge1993,Bish2004} for some examples in the airline sector, where the largest aircraft is assigned to the legs that expect the highest and most valuable demand.
	Interesting related contributions that explicitly address the capacity uncertainty can be found in  \cite{Busing2018,wang2006}. \color{black}Further recent applications that are somehow related to such setting are worth mentioning, that are mainly centered on pricing problems. We particularly refer the reader to the attended home delivery management \cite{klein18, klein19,mackert19,Yang16, Yang17} where the impact of routing costs has to be estimated in advance in order to evaluate the contribution margin of each request. \color{black}

	\section{Fundamentals of revenue management capacity-control}
	\label{section2}
	Multiple-resource capacity control problems arise when each product is sold as a set of resources.
	Denoting by ${\cal L}$ and ${\cal N}$ the index-sets of resources and products, respectively, with $|{\cal L}|=l$ and $|{\cal N}|=n$, products are represented by the columns of an $L \times N$ binary matrix  $\bm{\Gamma}$, whose component $\gamma_{ij}$ is set to one if resource $i \in {\cal L}$ is used by product $j\in {\cal N}$, and set to zero otherwise. Resource availability at the beginning of the time horizon is denoted by an $L$-dimensional vector $\bm{\xi}$, while $\mathbf{p}$ denotes the $n$-dimensional vector of the prices per unit-product.
	
	\color{black}Assuming that \color{black}the booking horizon is divided in $T$ time periods indexed by $t \in {\cal T} \triangleq \{1,\ldots,T\}$, we adopt as a standing assumption that at most one request for a product can arrive within each booking-period. The demand process is represented by a sequence of vectors $\{\mathbf{P}(t):t\geq 1\}$, where $\mathbf{P}(t)=\big(P_1(t),\ldots,P_n(t)\big)^\top$, and $P_1(t),\ldots,P_n(t)$ are random variables, independent across time and with known joint distribution, representing the demand of product $j$ at period $t$. In particular, for any product $j$ the demand random component $P_j(t)$ has two possible realizations, namely, $p_j$ and $0$. The case $P_j(t)=p_j$ means that a request for product $j$, with sell price equal to $p_j$, has been issued at period $t$, while $P_j(t)=0$ means that no request of product $j$ has been issued at period $t$. In case $P_j(t)=0$ for every $j$, the obvious meaning is that no product request for any product has been issued at period $t$. For any given pair $(t,j)$, we denote by $\lambda_j^t$ the probability that one request for product $j$ is issued at period $t$, and by $\lambda_0^t$ the probability that no request arrives at period $t$, hence $\lambda_0^t+\sum_{j \in {\cal N}} \lambda_j^t=1$.
	
	The decision variables of the problem are represented by a set of $n$-dimensional binary vectors $\{\mathbf{u}(t): t \in {\cal T}\}$, where $u_{j}(t)$  is set to one if and only if a request for product $j$ is accepted at period $t$. We adopt, as a system-state representation at any period $t$,  the vector $\mathbf{w}(t)=\big(w_1(t),\ldots,w_n(t)\big)^\top$ denoting the cumulative decisions made until the end of period $t-1$, namely, $w_j(t)$ represents the number of requests for product $j$  accepted before entering period $t$. As a consequence, if a request for product $j$ is accepted at period $t$, then the state $\mathbf{w}(t)$ is updated according to $\mathbf{w}(t+1)=\mathbf{w}(t)+\mathbf{e}_j$, where $\mathbf{e}_j$ is the $j$th unit vector.
	
	Observing that at any period $t$ the accept/reject decision depends on the state $\mathbf
	w(t)$ and on the prices $\mathbf{p}$, that is, $\mathbf{u}(t)=\mathbf{u}\big(t,\mathbf{w}(t),\mathbf{p}\big)$, the revenue management decision, in the framework of capacity control, is the following: at period $t$, given the current system-state $\mathbf{w}(t)$ and the random demand vectors $\mathbf{P}(t),\ldots,\mathbf{P}(T)$, should we accept or reject the incoming request?
	
	The problem can be formulated as a dynamic stochastic program, whose aim is to determine the optimal decisions $\mathbf{u}^*(t)=\mathbf{u}^*\big(t,\mathbf{w}(t),\mathbf{p}\big)$, where
	\begin{equation*}
	\mathbf{w}(t) \triangleq
	\begin{cases}
	\sum_{\tau = 1}^{t-1} \mathbf{u}^*(\tau), \,\, \mathrm{if} \,\, t \geq 2 \\[10pt]
	0, \,\, \mathrm{if} \,\, t = 1.
	\end{cases}
	\label{systemstate}
	\end{equation*}
	Denoting by $\rho_t(\mathbf{w})$ the maximum expected revenue\footnote{For simplicity of notation, the dependence on $t$ is dropped whenever possible.} from period $t$, with $\mathbf{w}$ representing the cumulative accepted requests until $t-1$, then $\rho_t(\mathbf{w})$ must satisfy the following Bellman equation, see \cite{bertsimas03}:
	\begin{equation}
	\rho_t(\mathbf{w}) = \lambda_0^t \rho_{t+1}(\mathbf{w}) + \sum_{j \in {\cal N}}\lambda_j^t \max_{u_j \in \{0,1\}}\Big\{p_j u_j+ \rho_{t+1}(\mathbf{w}+u_j \mathbf{e}_j ) \Big\},
	\label{bellman1}
	\end{equation}
	with the boundary conditions
	\begin{equation}
	\rho_{T+1} (\mathbf{w}) = 0 \;\; \forall \mathbf{w}: \bm{\xi} \geq \bm\Gamma \mathbf{w}
	\label{boundary1}
	\end{equation}
	and
	\begin{equation}
	\rho_{t}(\mathbf{w}) = -\infty \;\;\; \forall t ,\;\forall \mathbf{w}:  \bm{\xi} \not \geq \bm\Gamma \mathbf{w}.
	\label{boundary2}
	\end{equation}
	By observing that the Bellman equation (\ref{bellman1}) can be equivalently reformulated as
	\begin{equation}
	\rho_t(\mathbf{w}) = \rho_{t+1}(\mathbf{w}) + \sum_{j \in {\cal N}}\lambda_j^t \max\Big\{p_j - \rho_{t+1}(\mathbf{w}) + \rho_{t+1}(\mathbf{w}+\mathbf{e}_j),0 \Big\},
	\label{bellman2}
	\end{equation}
	the optimal control $\mathbf{u}^*(t)$ can be obtained according  to the criterion
	\begin{equation}
	u_j^*(t)=%u_j^*(t,\mathbf{w},p_j) =
	\begin{cases}
	1 & \mbox{if} \;\; p_j \geq \rho_{t+1}(\mathbf{w}) - \rho_{t+1}(\mathbf{w}+\mathbf{e}_j) \,\, \mathrm{and} \,\, \bm{\xi} \geq \bm\Gamma (\mathbf{w} + \mathbf{e}_j),\\[10pt]
	0 & \mbox{otherwise},
	\end{cases}
	\label{revenue-u-opt}
	\end{equation}
	namely, a request for product $j$ at time $t$ is accepted if and only if there is enough residual capacity, and the future maximum expected revenue if the request is accepted, $p_j + \rho_{t+1}(\mathbf{w}+\mathbf{e}_j)$, is not lower than the future maximum expected revenue if the request is rejected, $\rho_{t+1}(\mathbf{w})$.
	%\begin{remark}
	We note that it is easy to show that for any given $t$ the revenue $\rho_t(\cdot)$ is monotonically nonincreasing.
	
	%	Furthermore, we recall that the difference
	%	\begin{equation}
	%	\omega_j(\mathbf{w},t) \triangleq  \rho_{t+1}(\mathbf{w}) - \rho_{t+1}(\mathbf{w}+ \mathbf{e}_j) \geq 0 \label{opportunitycostnew}
	%	\end{equation}
	%	is usually referred to as the ``opportunity cost'' of product $j$ at time $t$, as it represents the cost (in terms of future expected-revenue loss) of accepting the request of product $j$ at time $t$. Hence, the criterion (\ref{revenue-u-opt}) claims that a request for product $j$ at time $t$ is accepted if and only if there is enough residual capacity, i.e., $\bm{\xi} \geq \bm\Gamma (\mathbf{w} + \mathbf{e}_j)$, and the corresponding price $p_j$ is not lower than its opportunity cost $\omega_j(\mathbf{w},t)$. %\Halmos
	%%\end{remark}
	
	The large size of the state space for meaningful real-life problems makes Bellman equations such as (\ref{bellman1}) impossible to solve and, as a consequence, optimal controls (\ref{revenue-u-opt}) impossible to calculate. For such reason, several approximate models have been introduced in the literature, based either on simplified formulations or on decompositions, along with several control policies. For later convenience, next we briefly review the Deterministic Linear Programming (DLP) model and the related control policies. Denoting by $D_j(t)$ the random demand of each product $j$ from period $t$ to period $T$, with expectation $\mathbb{E}[D_j(t)]$, the DLP approximation is formulated as
	\begin{equation}
	\widetilde{\rho}_t(\mathbf{w}) = \max_{\mathbf{y}} \left\{\mathbf{p}^\top \mathbf{y} : \bm\Gamma \mathbf{y} \leq \bm{\xi} - \bm\Gamma \mathbf{w}, \,  \mathbf{0} \leq \mathbf{y} \leq \bm{\mu} \right\},
	\label{LP-approx}
	\end{equation}
	where the demand appears as a vector of deterministic parameters
	$$
	\bm{\mu}\triangleq\big(\mathbb{E}[D_1(t)],\ldots,\mathbb{E}[D_n(t)]\big)^\top,
	$$
	and the decision variables $\mathbf{y} = \left(y_1,\ldots,y_n\right)$ represent the so called \textit{booking limit} of each product. Therefore, with respect to the current state $\mathbf{w}$ at time $t$, the optimal solution of problem (\ref{LP-approx}), say $\widetilde{\mathbf{y}}=\left(\widetilde{y}_1,\ldots,\widetilde{y}_n\right)$, can be seen as the most desirable realizations (with respect to future revenues) of the future-demand of each product, whose acceptance can be planned in advance. 

	\section{Capacity control for  distribution logistics problems} \label{section3} We aim at extending model (\ref{bellman1})-(\ref{boundary2}) to capacity control problems for distribution logistics. As announced earlier, we consider a scenario where a logistic service provider, along the given booking horizon ${\cal T}$, dynamically receives requests of parcel collection from a set of customers located at the nodes of a distribution network. In our setting we think of nodes as zones of the network where demand is clustered (e.g., districts in an urban network). Although the requests refer to a future operational time-horizon, still the accept/reject decision must be made at the time each request is received. Upon termination of the booking horizon, the logistic operator will optimally plan the collection operations, by solving an appropriate cost-minimization problem. The profit, associated to the whole process, will then be obtained by subtracting the operational optimal cost from the total revenues returned by the accepted requests.
	%, considering the operational feasibility and the expected future revenue associated to the resulting routing solution.
	
	The relevant updates with respect to the capacity-control setting reviewed in \S \ref{section2} are the following:
	\begin{itemize}
		\item The resources used in the process (vehicles, depots, drivers, etc.) depend on the specific structure of the distribution problem; %, and are typically represented by a fleet of vehicles (and possibly a crew of drivers);
		\item Products are associated to servicing customers; %i.e., visiting nodes with vehicles whose residual available space can accommodate the requested delivery;
		\item Products are no longer represented as a given set of resources, since the logistic operator has several options to arrange the collection plan (i.e., to consume resources), even choosing the minimal-cost one;
		\item The boundary conditions must appropriately account for the revenue loss associated to the operational cost of the service, and to possible infeasibility of the set of accepted requests.
	\end{itemize}
	We are now ready to set some notation. We assume that a complete directed graph ${\cal G} = ({\cal V},{\cal A})$ is given, where ${\cal V} = \{0,\ldots,n\}$ is the set of nodes and ${\cal A} \subseteq \{(i,j): i \in {\cal V},j \in {\cal V},i \neq j\}$ is the set of arcs. Each node in the set ${\cal N} = \{1,\ldots,n\}$ is associated to \color{black}a geographical zone shared by a set of customers, \color{black}while node $0$ corresponds to the depot, hence ${\cal V} = {\cal N} \cup \{0\}$. A nonnegative traveling cost $c_{ij}$ is associated with each arc $(i,j)\in {\cal A}$. A fleet of vehicles ${\cal K} = \{1,\ldots,K\}$ is available for visiting \color{black}nodes \color{black}in the network, each vehicle $k \in {\cal K}$ being associated to a capacity $Q_k$.  We assume, for simplicity of notation, and without loss of generality, that there is only one kind of item to be collected on the network, and that the number of items to be collected for each request is just one.
	% It is worth highlighting  that the correspondence between nodes and customers simplifies the product representation and the whole notation, since there is one product for each node-customer in the network, but it is not restrictive.
	%Let $\bm{\Gamma}$ be an $K \times n$ binary matrix, whose component $\gamma_{kj}$ is set to one if the vehicle $k \in {\cal K}$ is used to collect items from customer $j\in {\cal N}$, and set to zero otherwise. Let $\bm{\Gamma}_j$ be the $j$th column of $\bm{\Gamma}$, where $\bm{\Gamma}_j$ denotes the set of vehicles used to satisfy the requests from customer $j$.
	%Let $\bm{\xi}$ denote a $K$-dimensional vector containing the resource (i.e., the vehicle capacity) availability  at the beginning of the time horizon, that is $\bm{\xi}_k=Q_k$, $\forall k \in {\cal K}$.
	As a consequence, each demand random-variable $P_j(t)$ represents the demand of node $j$ at period $t$, whose collection price is denoted by $p_j$. As already stated in \S \ref{section2}, for each node $j$, the demand random component $P_j(t)$ has two possible realizations, namely, $p_j$ and $0$. The case $P_j(t)=p_j$ means that a one-item collection request at node $j$ with sell price equal to $p_j$ has been issued at period $t$, while $P_j(t)=0$ means that no request at node $j$ has been issued at period $t$. In case $P_j(t)=0$ for every $j$, the obvious meaning is that no item collection request from any customer has been issued at period $t$. Moreover, as previously stated, for any given pair $(t,j)$, we denote by $\lambda_j^t$ the probability that a one-item collection request at node $j$ is issued at period $t$, and by $\lambda_0^t$ the probability that no request arrives at period $t$, hence $\lambda_0^t+\sum_{j \in {\cal N}} \lambda_j^t=1$.

	The decision variables of the problem are again represented by a set of $n$-dimensional binary vectors $\{\mathbf{u}(t): t \in {\cal T}\}$, where $u_{j}(t)=1$ if and only if a one-item collection request at node $j$ is accepted at period $t$, and $u_{j}(t)=0$ otherwise. The system state, at any period $t \in {\cal T}$, is represented as the vector $\mathbf{w}(t)=\big(w_1(t),\ldots,w_n(t)\big)^\top$ denoting the cumulative decisions made until the end of period $t-1$.
	Indeed, the component $w_{j}(t) = \sum_{\tau = 1}^{t-1} u_{j}^*(\tau)$ represents the number of requests (i.e., the number of items to be collected) accepted until the end of the booking period $t-1$ at node $j$.
	%Recalling the definition of the matrix  $\bm{\Gamma}$,
	%Let $\bm{\xi}(1)$ represent the resource (i.e., total vehicle capacity) availability  at the beginning of the time horizon, that is $\bm{\xi}(1)=\sum_{k}^{K}Q_k$, it is easy to verify that the residual capacity at period $t$ can be calculated as $\bm\xi(t) = \bm\xi(1)-e^T \mathbf{w}(t)$.
	
	For any given system-state $\mathbf{w}$, the operational cost $v(\mathbf{w})$ associated with servicing the accepted requests $(w_1,\ldots,w_n)$ depends on the specific features of the underlying distribution problem. In an abstract setting, denoting by $\mathbf{x} \in X\big(\mathbf{w}\big)$ an $m$-dimensional vector representing the operational decision variables, by $X\big(\mathbf{w}\big) \subset \mathbb{R}^m$ the set of feasible operational decisions related to $\mathbf{w}$, and by $z : \mathbb{R}^m \mapsto \mathbb{R}$ the operational cost-function, $v(\mathbf{w})$ can be obtained as
	\begin{equation}
	v\big(\mathbf{w}\big) \triangleq
	\begin{cases}
	z^*(\mathbf{w}), \mbox{ if } X\big(\mathbf{w}\big) \neq \emptyset \\
	+ \infty, \mbox{ if } X\big(\mathbf{w}\big) = \emptyset,
	\end{cases}
	\label{subproblem1:static}
	\end{equation}
	where
	\begin{equation}
	z^*(\mathbf{w}) = \min \Big\{z(\mathbf{x}): \mathbf{x} \in X\big(\mathbf{w}\big) \Big\}.
	\label{subproblem2:static}
	\end{equation}
	In the following, we will work under the assumptions that there exists a finite minimum $z^*\big(\mathbf{w}\big) \geq 0$ if $X\big(\mathbf{w}\big) \neq \emptyset$, that $X\big(\mathbf{0}\big) \neq \emptyset$, and that $z^*(\mathbf{0})=0$.
	
	By letting $\pi_t(\cdot)$ the maximum expected profit from booking period $t$, with respect to the system state $\mathbf{w}$, the capacity control problem can be formulated as the following dynamic program:
	\begin{equation}
	\pi_t(\mathbf{w}) = \lambda_0^t \pi_{t+1}(\mathbf{w}) + \sum_{j \in {\cal N}}\lambda_j^t \max_{u_j \in \{0,1\}}\big\{p_j u_j+ \pi_{t+1}(\mathbf{w}+u_j \mathbf{e}_j ) \big\},
	\label{bellman1new}
	\end{equation}
	with the boundary condition
	\begin{equation}
	\pi_{T+1}(\mathbf{w})= -v(\mathbf{w})
	\label{boundary1new}
	\end{equation}
	for which an optimal control $\mathbf{u}^*(t)$ can be obtained according  to
	\begin{equation}
	u_j^*(t)=%u_j^*(t,\mathbf{w},p_j) =
	\begin{cases}
	1 & \mbox{if} \;\; p_j \geq \pi_{t+1}(\mathbf{w}) - \pi_{t+1}(\mathbf{w}+\mathbf{e}_j) \,\, \mbox{and} \,\, X\big(\mathbf{w}\big) \neq \emptyset \\[10pt]
	0 & \mbox{otherwise}.
	\end{cases}
	\label{revenue-u-opt-new}
	\end{equation}
	
	For a deeper understanding of the dynamic program (\ref{bellman1new})-(\ref{boundary1new}) and to explore the practical role of such framework, now we adopt a CVRP as a specific distribution application. In particular, we assume that the logistics service provider, for a given system-state $\mathbf{w}$, has to plan the routes of at most $K$ vehicles with given capacity, in order to collect $w_j$ items at each node $j$ on the network. In the following, we provide the details of the adopted CVRP formulation. Denoting by $\alpha_{ij}^k$ a binary decision variable set to one if arc $(i,j)$ is used by vehicle $k$, by $\beta_{i}^k$ a binary decision variable set to one if node $i \in {\cal V}$ is visited by vehicle $k$, and by $q_{j}^k$ a nonnegative decision variable representing the number of items collected at node $j \in {\cal N}$ by vehicle $k$, we formulate the operational problem as follows:
	
	\begingroup
	\addtolength{\jot}{1.1em}
	\begin{alignat}{3}
	z^*(\mathbf{w}) = \text{minimize} {} \qquad  & \sum_{k \in {\cal K}} \sum_{(i,j) \in {\cal A}} c_{ij} \alpha_{ij}^k  &  \label{routing-cost} \\
	\text{subject to} \qquad & \sum_{j:(i,j) \in {\cal A}} \alpha_{ij}^{k} = \beta_i^k \qquad & \forall i \in {\cal V}, \, \forall k \in {\cal K} \label{outflow}\\
	\qquad & \sum_{j:(j,i) \in {\cal A}} \alpha_{ji}^{k} = \beta_i^k \qquad & \forall i \in {\cal V}, \, \forall k \in {\cal K} \label{inflow}\\
	\qquad & \sum_{k \in {\cal K}} \beta_{j}^{k} \leq  1 \qquad & \forall j \in {\cal N} \label{assignment}\\
	\qquad & \sum_{k \in {\cal K}} \beta_{0}^{k} \leq K & \label{fleet} \\
	\qquad & \sum_{i \in {\cal S}} \sum_{j \in {\cal V}\setminus{\cal S}} \alpha_{ij}^k \geq \beta_h^k \qquad & \forall {\cal S} \subset {\cal V}:0 \in {\cal S}, \forall h \in {\cal V}\setminus {\cal S}, \forall k \in {\cal K} \label{subtour}\\
	\qquad & q_j^k \leq Q_k \beta_j^k \qquad & \forall j \in {\cal N}, \, \forall k \in {\cal K} \label{vehiclecapacity}\\
	\qquad & \sum_{j \in {\cal N}} q_j^k \leq Q_k \qquad & \forall k \in {\cal K} \label{routecapacity}\\
	\qquad & \sum_{k \in {\cal K}} q_j^k = w_j \qquad & \forall j \in {\cal N} \label{requests}\\
	\qquad & \alpha_{ij}^k \in \{0,1\} \qquad & \forall (i,j) \in {\cal A}, \forall k \in {\cal K} \label{flowvars}\\
	\qquad & \beta_{i}^k \in \{0,1\} \qquad & \forall i \in {\cal V}, \forall k \in {\cal K} \label{assignvars}\\
	\qquad & q_{j}^k \geq 0 \qquad & \forall j \in {\cal N}, \forall k \in {\cal K} \label{amountvars}
	\end{alignat}
	The objective function (\ref{routing-cost}) minimizes the total routing cost, while constraints (\ref{outflow}) and (\ref{inflow}) guarantee that an arc enters and an arc leaves from each node which is visited. Constraints (\ref{assignment}) ensure that each node is visited at most once, constraint (\ref{fleet}) limits the number of routes (i.e., vehicles) to be at most $K$, and constraints (\ref{subtour}) impose that the solution is connected, with constraints (\ref{flowvars}), (\ref{assignvars}), and (\ref{amountvars}) representing variables definitions. The side constraints (\ref{vehiclecapacity}) to (\ref{requests}) concern the fulfillment of the requests. In particular, constraints (\ref{vehiclecapacity}) ensure that when a vehicle visits a node it cannot collect an amount larger than its capacity, constraints (\ref{routecapacity}) ensure that the total amount of items collected by each vehicle is not larger than its capacity, and constraints (\ref{requests}) ensure that the requests at each node are fulfilled.
	%In Table \ref{vrpnotation} we summarize the notation for all data and variables of problem (\ref{routing-cost})--(\ref{amountvars}).
	Finally, we observe that by letting $\mathbf{x} = (\bm{\alpha},\bm{\beta},\mathbf{q})$ and
	$$
	X(\mathbf{w})=\{\mathbf{x} : \mathbf{x}=(\bm{\alpha},\bm{\beta},\mathbf{q}) \,\, \mathrm{and} \,\, (\ref{outflow})-(\ref{amountvars}) \,\, \mathrm{hold} \}
	$$
	then the CVRP (\ref{routing-cost})--(\ref{amountvars}) has the same structure as problem  (\ref{subproblem2:static}), and that the existence of a finite minimum $z^*(\mathbf{w}) \geq 0$ if $X(\mathbf{w}) \neq \emptyset$ is ensured.

	\section{Deterministic approximation and booking limit control}\label{section4}
	The proposed dynamic program (\ref{bellman1new})-(\ref{boundary1new}) cannot be solved optimally due to the curse of dimensionality. Hence, similar to what reviewed in \S \ref{section2}, we introduce a deterministic approximation of $\pi_t(\cdot)$, whose solution provides the key information for developing approximate control policies, that support the logistic operator in the real-time decision-making process. In fact, although the deterministic approximation is a static model, still it can be solved dynamically along the booking horizon by appropriately updating both demand and system-state information.
	
	First we focus on the abstract setting of the dynamic program (\ref{bellman1new})-(\ref{boundary1new}). For any given booking period $t$, we recall that
	the requests accepted prior to $t$ cannot be turned down in the future, hence the future profit depends on the random revenues associated to the requests accepted in the future, and on the optimal distribution cost of the whole set of accepted requests. Denoting by
	$$
	\bm{\mu}(t)=\big(\mathbb{E}[D_1(t)],\ldots,\mathbb{E}[D_n(t)]\big)^\top
	$$
	the vector of future demand expectations at each node $j$, and introducing for each $j$ a booking-limit decision variable $y_j$, we follow the guidelines introduced in \S \ref{section2} to formulate a deterministic approximation, $\widetilde{\pi}_t(\mathbf{w})$, with respect to a feasible system-state $\mathbf{w}$ as
	\begin{equation}
	\widetilde{\pi}_t(\mathbf{w})=\max_{\mathbf{x},\mathbf{y}} \big\{\mathbf{p}^\top \mathbf{y} - z(\mathbf{w}+\mathbf{y}): \mathbf{x} \in X(\mathbf{w}+\mathbf{y}), \,\, \mathbf{0} \leq \mathbf{y} \leq \bm{\mu}\big\}
	\tag{$\widetilde{P}_t\left(\mathbf{w}\right)$}
	\label{det-approx}
	\end{equation}
	and we denote by $\big(\widetilde{\mathbf{x}}(\mathbf{w}),\widetilde{\mathbf{y}}(\mathbf{w})\big)$ the optimal solution of \ref{det-approx}. Therefore, with respect to the current state $\mathbf{w}$ at booking period $t$, the optimal solution vector $\widetilde{\mathbf{y}}$ can be seen as the most desirable combination (with respect to future profits) of the random future-demand realizations of each product.
	
	The practical role of deterministic approximations is that they provide the relevant information needed to build approximate control policies as  the \textit{partitioned booking limits} described next.
	
	A booking-limit control policy consists in setting to $\widetilde{y}_j$ the upper bound on the number of requests of each product $j$ that can be accepted from the current period $t$ until $T$. It corresponds then to allocating a fixed amount of capacity on the involved resources.
	%For further details on the role of the deterministic approximation (\ref{LP-approx}), along with its theoretical properties and practical applications, we refer the reader to \cite{talluri05} and the references therein.
	We propose in Algorithm \ref{bookinglimit} a detailed procedure for the dynamic accept/reject decision process based on booking-limit controls. The procedure is based on selecting a subset ${\cal S} \subseteq {\cal T}$ of booking periods (referred to as solution periods) when the deterministic approximation \ref{det-approx} is actually solved. By iterating over all booking periods, whenever a request arrives at $t$, the accept/reject decision is based on the residual booking limit $\overline{\mathbf{y}}$ related to the most recent solution of \ref{det-approx}, see steps \ref{Step1:7}-\ref{Step1:12}. At the end of the booking horizon, the total profit is obtained, at step \ref{Step1:17}, by preliminary solving the operational cost-minimization problem at step \ref{Step1:16}.
	\begin{algorithm}
		\scriptsize
		\caption{Booking Limit Policy (BL)}\label{bookinglimit}
		\begin{algorithmic}[1]
			\algrenewcommand\algorithmicrequire{\textbf{Input:}}
			\algrenewcommand\algorithmicensure{\textbf{Output:}}
			\Require the set of booking periods ${\cal T} = \{1,\ldots,T\}$ and the set of solution periods ${\cal S} \subseteq  {\cal T}$, with $\{1\} \subseteq {\cal S}$
			\Ensure a system state $\mathbf{w}$ along with the corresponding profit $\pi^*$
			%\Statex
			\State{Set $\mathbf{w}=\mathbf{0}$ and $t = 1$} \label{Step1:1} \Comment{Initialization}
			\If{$t\in {\cal S}$} \label{Step1:2} \Comment{$t$ is a solution period}
			\State set $\bm{\mu}(t)=\mathbb{E}[\mathbf{D}(t)]$ \Comment{Update the future mean-demand vector} \label{Step1:3}
			\State solve \ref{det-approx} and obtain $\widetilde{\mathbf{y}}(\mathbf{w})$ \label{Step1:4} \Comment{Calculate the booking limit}
			\State set $\overline{\mathbf{y}}=\widetilde{\mathbf{y}}(\mathbf{w})$ \label{Step1:5} \Comment{Store the residual booking limit}
			\EndIf \label{Step1:6}
			\If{product $j$ is requested at time $t$} \label{Step1:7}
			\If{$\overline{y}_j \geq 1$} \label{Step1:8} \Comment{Check feasibility of booking limit}
			\State set $\mathbf{w} = \mathbf{w} + \mathbf{e}_j$ \label{Step1:9} \Comment{Update the system state}
			\State set $\overline{\mathbf{y}} = \overline{\mathbf{y}} -\mathbf{e}_j$ \label{Step1:10} \Comment{Update the residual booking limit}
			\EndIf \label{Step1:11}
			\EndIf \label{Step1:12}
			\If{$t < T$} \label{Step1:13}
			\State set $t = t+1$ and \textbf{go to \ref{Step1:2}} \label{Step1:14}
			\Else  \Comment{End of booking horizon reached} \label{Step1:15}
			\State calculate $z^*(\mathbf{w})$ by solving (\ref{subproblem2:static}) \label{Step1:16} \Comment{Calculate the operational cost}
			\State set $\pi^* = \mathbf{p}^\top \mathbf{w} - z^*(\mathbf{w})$ and \textbf{exit} \label{Step1:17} \Comment{Calculate the total profit}
			\EndIf \label{Step1:18}
		\end{algorithmic}
	\end{algorithm}

	A specific version of the deterministic approximation \ref{det-approx}, and of the booking limit policy described above, can be obtained by focusing  on the CVRP-based capacity control dynamic-program (\ref{bellman1new})-(\ref{boundary1new}), where we  adopt (\ref{routing-cost})-(\ref{amountvars}) as the cost-minimization problem in (\ref{subproblem2:static}). The resulting specific formulation of \ref{det-approx} is the following \textit{Profit Maximization Vehicle Routing Problem} (PMVRP):
	
	\begingroup
	\addtolength{\jot}{1.1em}
	\begin{alignat}{3}
	\widetilde{\pi}_t(\mathbf{w}) = \text{maximize} {} \qquad  & \sum_{j \in {\cal N}} p_j y_j - \sum_{k \in {\cal K}} \sum_{(i,j) \in {\cal A}} c_{ij} \alpha_{ij}^k  &  \label{RMVRP:obj} \\
	\text{subject to} \qquad & (\ref{outflow})-(\ref{routecapacity}), (\ref{flowvars})-(\ref{amountvars}) \nonumber \\
	\qquad & \sum_{k \in {\cal K}} q_j^k = w_j + y_j \qquad & \forall j \in {\cal N} \label{newrequests}\\
	\qquad & 0 \leq y_j \leq \mu_j \qquad & \forall j \in {\cal N} \label{meandemand}
	\end{alignat}
	a mixed-integer linear program, where constraints (\ref{newrequests}) ensure that at least the requests accepted prior to $t$ are fulfilled.
	%The solution of the PMVRP allows to apply policies like booking limits and bid-price controls following the steps introduced in Algorithms 1 and 2.
	We remark that unlike common VRP formulations that are typically focused either on cost minimization or on prize maximization (the so called VRP with profits), see \cite{toth14}, the objective function of PMVRP is the difference between the revenues associated to the accepted requests and the routing costs associated to the fulfillment of such requests. Next, we observe that, focusing on the constraint set, the PMVRP formulation looks as a CVRP with side constraints, sharing with standard CVRP formulations the sets of constraints  (\ref{outflow})-(\ref{routecapacity}) and (\ref{flowvars})-(\ref{amountvars}). We also note that constraints (\ref{requests}) are replaced by (\ref{newrequests}) and (\ref{meandemand}), both representing the relevant features of a RM formulation, where $y_j$ is the most desirable demand realization, not larger than the future demand expectation, at each node $j$. We finally notice that, to the best of our knowledge, the PMVRP formulation appears as a new model in the field of vehicle routing problems that would certainly require an in-depth methodological analysis.

	\section{Computational experience}
	\label{section5}
	Our goal is to show the practical relevance of the RM technique introduced in \S \ref{section4}, by testing it on the specific CVRP-based application described above. Therefore, we will next refer to the deterministic approximation \ref{det-approx} as the PMVRP formulation (\ref{RMVRP:obj})-(\ref{meandemand}), and we will accordingly update the steps of the approximate control policy.
	
	The computational issues posed by the approximate control policy are essentially related to the need of repeatedly solving the PMVRP deterministic approximations during the booking process, and to solve a cost-minimization CVRP at the end of the booking horizon in order to calculate the total profit. This computational burden can be significantly softened in Algorithm \ref{bookinglimit} by keeping the cardinality of ${\cal S}$ sufficiently small. \color{black}As a consequence, in our implementation of Algorithm \ref{bookinglimit} we have set ${\cal S}= \{1\}$ and ${\cal S}= \{1,2\}$. In the former case, next referred to as BLP, we solve PMVRP only once at the beginning of the booking horizon. In the latter case, next referred to as BLPR (i.e., booking limit policy with reoptimization), PMVRP is solved twice, at the beginning and in the middle of the booking horizon. Thus, by adopting an exact solver equipped with some truncation criterion, we manage to obtain results in a reasonable amount of time.

%	 For a comparison against the proposed techniques, we will consider two alternative approaches, one based on the myopic First-Come First-Served strategy, the other on a Perfect Knowledge assumption. Before describing several implementation details of each tested method and the related results, we first introduce some relevant information about the generation of the adopted instance test-set.

	\subsection{Test examples}
%	We have first considered two toy examples with $n=5$ and $n=10$ nodes, a number of vehicles $K=3$ and $5$ respectively and a capacity $Q_k=20$ for each vehicle $k\in K$. We named the two examples as $T5$ and $T10$. These examples are useful since they can be used to test our policies by solving to   optimality, in a reasonable amount of time, both NRMVRP and CVRP.
We have generated three groups of instances, related to the cases $n \in \{15, 25, 50\}$. For each group we have generated 3 instances, by selecting the network topologies, and the related initial expected-demand $\mu_j$ at each node $j$, from the well known \textit{Random} ($R$),\textit{Clustered} ($C$) and \textit{Random}-\textit{Clustered} ($RC$) instances introduced in \cite{solomon1987}. The price $p_j$, at each node $j$, has been set according to an inverse proportionality rule with respect to the expected demand. For every instance the number of  vehicles $K$, and the capacity $Q_k=Q$ for each $k$, have been selected so that the load factor 
$$
LF = \frac{\sum_{j=1}^n \mu_j}{kQ}
$$ 
is between $1.0$ and $1.5$. We refer to each instance adopting the notation $I.n$, where $I \in \{C,R,RC\}$ and $n \in \{15, 25, 50\}$. 
		
For each instance, we have generated 50 request paths in order to simulate the booking process 50 times. Each request path has been randomly generated in two steps. The first one aims at randomly sampling the number of requests $D_j(1)$, at each node $j$, from a truncated normal distribution with mean $\mu_j$ and coefficient of variation 0.1.  The second one aims at randomly sampling booking-request arrival-times from a uniform distribution between 1 and $T$.

	\subsection{Computational results}
	For a comparison against our two implementations of Algorithm \ref{bookinglimit}, i.e. BLP and BLPR, we will consider two alternative approaches, one based on the myopic First-Come First-Served (FCFS) strategy, the other on a Perfect Knowledge (PK) assumption. In particular, the FCFS strategy requires for every request to certify the existence of a feasible route fulfilling all the past accepted requests plus the new one, and only in such case to accept the incoming request. Moreover, at the end of the booking horizon, in order to calculate the total profit, it will always be necessary to solve the CVRP (\ref{routing-cost})-(\ref{amountvars}). Unlike the myopic FCFS,  the PK approach borrows much of structure of the booking-limit control policy. In particular, it consists of an implementation of the \textit{booking-limit} policy, described in Algorithm \ref{bookinglimit}, by substituting the mean $\mu_j(t)$ by the actual realization $D_j(t)$ made dynamically available during the simulation process.
	All numerical experiments have been carried out on a Intel(R) Core(TM)i7-4510U CPU 2.60 GHz processor, implementing the four control policies BLP, BLPR, FCFS, and PK, in MATLAB r2017b. \color{black}We have adopted ILOG CPLEX 12.8 as a mixed-integer linear programming solver, and we have tuned the truncation criterion, based on the number of incumbent updates, by running few experiments on very small-size instances that can also be solved at optimality. As a result we decided to truncate the MIP solver after two incumbent updates when solving PMVRP, and after 4 incumbent updates when solving CVRP.
	
	A summary of the results obtained from the execution of the four methods is available from Figure \ref{C15} to Figure \ref{RC50}, where for each instance we report the empirical cumulative distribution function (ECDF) of the profit with respect to the 50 simulation runs. 
	
	Generally speaking we observe that the results returned by the BLP approach are the best in terms of robustness, as shown by the small variance obtained in every instance. Focusing now on the comparison between BLP and BLPR, we see that on almost all cases BLP performs better than BLPR in terms of profit. The same observation holds when comparing BLP against FCFS. We recall that PKP can be expected to return a kind of upper bound on the average profit, given that it exploits the actual realizations of the demand, not its mean value.	The comparison of BLP against PKP shows that, as expected, PKP performs clearly better that BLP in terms of solution quality, although the former returns the worst results in terms of robustness, as the very large variance shows in every instance. The convincing performance obtained by the BLP approach are confirmed when considering execution time. In fact, the average execution time is less than 15 seconds, 18 seconds and 160 seconds, respectively, for the instances with $n=15$, $n=25$, and $n=50$. Such results are slightly better than BLPR, and much better than both PKP and FCFS. 
	\color{black}
		
\begin{figure}
	\begin{tikzpicture}
	\begin{axis}[height= 8cm,width=12cm,xlabel=$\pi^*$, ylabel=fraction of runs with profit within $\pi^*$,ymin=0.0,ymax=1.0,xmin=18000,xmax=350000,
	legend style={at={(0.5,-0.15)}, anchor=north,legend columns=-1}
	]
	\addplot+ [
	%smooth,
	] coordinates {
		(300399.842, 0.000)
		(306951.092, 0.020)
		(310468.000, 0.100)
		(313984.908, 0.220)
		(317501.817, 0.400)
		(321018.725, 0.740)
		(324535.633, 1.000)
	};
	
	\addplot+ [
	%smooth,
	] coordinates {
		(235325.544, 0.000)
		(252711.575, 0.200)
		(267720.581, 0.480)
		(282729.587, 0.740)
		(297738.592, 0.900)
		(312747.598, 0.960)
		(327756.604, 1.000)
	};

	\addplot+ [
	%smooth,
	] coordinates {
		(229837.532, 0.000)
		(242550.403, 0.100)
		(252941.683, 0.260)
		(263332.963, 0.520)
		(273724.243, 0.880)
		(284115.524, 0.960)
		(294506.804, 1.000)
	};
	
	\addplot+ [
	%smooth,
	] coordinates {
		(18843.287,	0.000)
		(73167.390,	0.120)
		(127301.157, 0.240)
		(181434.924, 0.340)
		(235568.691, 0.620)
		(289702.458, 0.640)
		(343836.225, 1.000)
	};

	\addlegendentry{BLP}
	\addlegendentry{BLPR}
	\addlegendentry{FCFS}
	\addlegendentry{PKP}
	\end{axis}
	\end{tikzpicture}
	\caption{Computational results: ECDF for instance $C.15$}
	\label{C15}
\end{figure}

\begin{figure}
	\begin{tikzpicture}
	\begin{axis}[height= 8cm,width=12cm,xlabel=$\pi^*$, ylabel=fraction of runs with profit within $\pi^*$,ymin=0.0,ymax=1.0,xmin=4400,xmax=60000,
	legend style={at={(0.5,-0.15)}, anchor=north,legend columns=-1}
	]
	\addplot+ [
	%smooth,
	] coordinates {
		(48995.563,	0.000)
		(49805.465,	0.040)
		(50120.462,	0.100)
		(50435.459,	0.120)
		(50750.456, 0.320)
		(51065.453, 0.500)
		(51380.450, 1.000)
	};
	
	\addplot+ [
	%smooth,
	] coordinates {
		(36678.726, 0.000)
		(40425.606, 0.060)
		(43801.994, 0.440)
		(47178.381, 0.640)
		(50554.769, 0.700)
		(53931.157, 0.900)
		(57307.544, 1.000)
	};
	
	\addplot+ [
	%smooth,
	] coordinates {
		(37209.487, 0.000)
		(39183.489, 0.160)
		(40781.637, 0.400)
		(42379.784, 0.640)
		(43977.932, 0.800)
		(45576.080, 0.960)
		(47174.228, 1.000)
	};
	
	\addplot+ [
	%smooth,
	] coordinates {
		(4453.244,	0.000)
		(13590.429,	0.020)
		(22682.632, 0.100)
		(31774.834, 0.240)
		(40867.037, 0.340)
		(49959.240, 0.400)
		(59051.443, 1.000)
	};
	
	\addlegendentry{BLP}
	\addlegendentry{BLPR}
	\addlegendentry{FCFS}
	\addlegendentry{PKP}
	\end{axis}
	\end{tikzpicture}
	\caption{Computational results: ECDF for instance $R.15$}
	\label{R15}
\end{figure}	
\begin{figure}
	\begin{tikzpicture}
	\begin{axis}[height= 8cm,width=12cm,xlabel=$\pi^*$, ylabel=fraction of runs with profit within $\pi^*$,ymin=0.0,ymax=1.0,xmin=64000,xmax=280000,
	legend style={at={(0.5,-0.15)}, anchor=north,legend columns=-1}
	]
	\addplot+ [
	%smooth,
	] coordinates {
		(207782.041,0.000)
		(212976.428, 0.040)
		(216072.006, 0.060)
		(219167.585, 0.240)
		(222263.163, 0.420)
		(225358.742, 0.700)
		(228454.320, 1.000)
	};

	\addplot+ [
	%smooth,
	] coordinates {
		(182740.066, 0.000)
		(191751.324, 0.020)
		(198916.722, 0.180)
		(206082.120, 0.500)
		(213247.518, 0.840)
		(220412.917, 0.980)
		(227578.315, 1.000)
	};

	\addplot+ [
	%smooth,
	] coordinates {
		(200321.438, 0.000)
		(207724.382, 0.102)
		(213103.877, 0.184)
		(218483.372, 0.551)
		(223862.868, 0.857)
		(229242.363, 0.939)
		(234621.859, 1.000)
	};

	\addplot+ [
	%smooth,
	] coordinates {
		(23116.582,	0.000)
		(64622.023,	0.020)
		(105893.962, 0.080)
		(147165.901, 0.160)
		(188437.840, 0.260)
		(229709.780, 0.320)
		(270981.719, 1.000)
	};

	\addlegendentry{BLP}
	\addlegendentry{BLPR}
	\addlegendentry{FCFS}
	\addlegendentry{PKP}
	\end{axis}
	\end{tikzpicture}
	\caption{Computational results: ECDF for instance $RC.15$}
	\label{RC15}
\end{figure}	

\begin{figure}
	\begin{tikzpicture}
	\begin{axis}[height= 8cm,width=12cm,xlabel=$\pi^*$, ylabel=fraction of runs with profit within $\pi^*$,ymin=0.0,ymax=1.0,xmin=150000,xmax=660000,
	legend style={at={(0.5,-0.15)}, anchor=north,legend columns=-1}
	]
	\addplot+ [
	%smooth,
	] coordinates {
		(538951.093, 0.000)
		(550456.061, 0.060)
		(556517.079, 0.120)
		(562578.097, 0.260)
		(568639.115, 0.560)
		(574700.132, 0.820)
		(580761.150, 1.000)
	};
	
	\addplot+ [
	%smooth,
	] coordinates {
		(483713.349, 0.000)
		(500073.138, 0.160)
		(511546.934, 0.400)
		(523020.730, 0.620)
		(534494.527, 0.760)
		(545968.323, 0.900)
		(557442.119, 1.000)
	};

	\addplot+ [
	%smooth,
	] coordinates {
		(392140.5327, 0.000)
		(414767.114, 0.020)
		(433432.681, 0.120)
		(452098.248, 0.360)
		(470763.815, 0.660)
		(489429.382, 0.920)
		(508094.949, 1.000)
	};
	
	\addplot+ [
	%smooth,
	] coordinates {
		(158519.115,	0.000)
		(242478.811,	0.040)
		(324837.305, 0.240)
		(407195.798, 0.520)
		(489554.292, 0.720)
		(571912.785, 0.860)
		(654271.279, 1.000)
	};

	\addlegendentry{BLP}
	\addlegendentry{BLPR}
	\addlegendentry{FCFS}
	\addlegendentry{PKP}
	\end{axis}
	\end{tikzpicture}
	\caption{Computational results: ECDF for instance $C.25$}
	\label{C25}
\end{figure}	

\begin{figure}
	\begin{tikzpicture}
	\begin{axis}[height= 8cm,width=12cm,xlabel=$\pi^*$, ylabel=fraction of runs with profit within $\pi^*$,ymin=0.0,ymax=1.0,xmin=2500,xmax=96000,
	legend style={at={(0.5,-0.15)}, anchor=north,legend columns=-1}
	]
	\addplot+ [
	%smooth,
	] coordinates {
		(65731.621,	0.000)
		(66827.287,	0.020)
		(67258.997, 0.120)
		(67690.708,	0.260)
		(68122.418, 0.460)
		(68554.128, 0.860)
		(68985.838, 1.000)
	};

	\addplot+ [
	%smooth,
	] coordinates {
		(55261.178, 0.000)
		(59906.250, 0.020)
		(63993.128, 0.080)
		(68080.006, 0.280)
		(72166.883, 0.800)
		(76253.761, 0.940)
		(80340.639, 1.000)
	};
	
	\addplot+ [
	%smooth,
	] coordinates {
		(56139.384, 0.000)
		(60139.588, 0.080)
		(63572.728, 0.240)
		(67005.868, 0.420)
		(70439.008, 0.660)
		(73872.147, 0.900)
		(77305.287, 1.000)
	};
	
	\addplot+ [
	%smooth,
	] coordinates {
		(2511.312, 0.000)
		(18110.073, 0.060)
		(33683.468, 0.100)
		(49256.863, 0.200)
		(64830.258, 0.520)
		(80403.653, 0.660)
		(95977.048, 1.000)
	};
	
	\addlegendentry{BLP}
	\addlegendentry{BLPR}
	\addlegendentry{FCFS}
	\addlegendentry{PKP}
	\end{axis}
	\end{tikzpicture}
	\caption{Computational results: ECDF for instance $R.25$}
	\label{R25}
\end{figure}		

\begin{figure}
	\begin{tikzpicture}
	\begin{axis}[height= 8cm,width=12cm,xlabel=$\pi^*$, ylabel=fraction of runs with profit within $\pi^*$,ymin=0.0,ymax=1.0,xmin=9600,xmax=500000,
	legend style={at={(0.5,-0.15)}, anchor=north,legend columns=-1}
	]
	\addplot+ [
	%smooth,
	] coordinates {
		(444316.173,	0.000)
		(452608.685,	0.040)
		(456413.154, 0.080)
		(460217.624, 0.220)
		(464022.093, 0.540)
		(467826.562, 0.840)
		(471631.032, 1.000)
	};
	
	\addplot+ [
	%smooth,
	] coordinates {
		(345276.810, 0.000)
		(361376.957, 0.040)
		(373989.459, 0.100)
		(386601.962, 0.240)
		(399214.464, 0.640)
		(411826.966, 0.940)
		(424439.468, 1.000)
	};
	
	\addplot+ [
	%smooth,
	] coordinates {
		(307618.344, 0.000)
		(319035.078, 0.080)
		(327344.556, 0.220)
		(335654.034, 0.480)
		(343963.512, 0.720)
		(352272.990, 0.920)
		(360582.468, 1.000)
	};
	
	\addplot+ [
	%smooth,
	] coordinates {
		(9601.881, 0.000)
		(91131.136, 0.020)
		(172563.402, 0.020)
		(253995.668, 0.060)
		(335427.934, 0.320)
		(416860.200, 0.460)
		(498292.466, 1.000)
	};

	\addlegendentry{BLP}
	\addlegendentry{BLPR}
	\addlegendentry{FCFS}
	\addlegendentry{PKP}
	\end{axis}
	\end{tikzpicture}
	\caption{Computational results: ECDF for instance $RC.25$}
	\label{RC25}
\end{figure}			
\begin{figure}
	\begin{tikzpicture}
	\begin{axis}[height= 8cm,width=12cm,xlabel=$\pi^*$, ylabel=fraction of runs with profit within $\pi^*$,ymin=0.0,ymax=1.0,xmin=165000,xmax=375000,
	legend style={at={(0.5,-0.15)}, anchor=north,legend columns=-1}
	]
	\addplot+ [
	%smooth,
	] coordinates {
		(298339.618, 0.000)
		(303599.868, 0.040)
		(305846.587, 0.100)
		(308093.305, 0.240)
		(310340.023, 0.460)
		(312586.742, 0.800)
		(314833.460, 1.000)
	};

	\addplot+ [
	%smooth,
	] coordinates {
		(245718.925, 0.000)
		(255284.975, 0.040)
		(262369.017, 0.140)
		(269453.058, 0.380)
		(276537.099, 0.680)
		(283621.141, 0.920)
		(290705.182, 1.000)
	};

	\addplot+ [
	%smooth,
	] coordinates {
		(265683.749, 0.000)
		(273344.525, 0.170)
		(278321.626, 0.383)
		(283298.728, 0.660)
		(288275.829, 0.872)
		(293252.931, 0.936)
		(298230.032, 1.000)
	};

	\addplot+ [
	%smooth,
	] coordinates {
		(162329.338, 0.000)
		(198963.670, 0.040)
		(233958.312, 0.120)
		(268952.953, 0.180)
		(303947.595, 0.260)
		(338942.236, 0.420)
		(373936.878, 1.000)
	};
	
	\addlegendentry{BLP}
	\addlegendentry{BLPR}
	\addlegendentry{FCFS}
	\addlegendentry{PKP}
	\end{axis}
	\end{tikzpicture}
	\caption{Computational results: ECDF for instance $C.50$}
	\label{C50}
\end{figure}

\begin{figure}
	\begin{tikzpicture}
	\begin{axis}[height= 8cm,width=12cm,xlabel=$\pi^*$, ylabel=fraction of runs with profit within $\pi^*$,ymin=0.0,ymax=1.0,xmin=220000,xmax=430000,
	legend style={at={(0.5,-0.15)}, anchor=north,legend columns=-1}
	]
	\addplot+ [
	%smooth,
	] coordinates {
		(360407.266, 0.000)
		(365755.711, 0.020)
		(367463.678, 0.040)
		(369171.645, 0.160)
		(370879.612, 0.400)
		(372587.579, 0.760)
		(374295.547, 1.000)
	};
	
	\addplot+ [
	%smooth,
	] coordinates {
		(275961.784, 0.000)
		(290725.032, 0.080)
		(302700.788, 0.240)
		(314676.544, 0.400)
		(326652.300, 0.620)
		(338628.056, 0.740)
		(350603.811, 1.000)
	};

	\addplot+ [
	%smooth,
	] coordinates {
		(251436.191, 0.000)
		(268408.398, 0.043)
		(282840.846, 0.085)
		(297273.294, 0.277)
		(311705.742, 0.574)
		(326138.190, 0.894)
		(340570.638, 1.000)
	};

	\addplot+ [
	%smooth,
	] coordinates {
		(229729.178, 0.000)
		(265002.799, 0.040)
		(297955.924, 0.060)
		(330909.049, 0.120)
		(363862.173, 0.260)
		(396815.298, 0.560)
		(429768.423, 1.000)
	};

	\addlegendentry{BLP}
	\addlegendentry{BLPR}
	\addlegendentry{FCFS}
	\addlegendentry{PKP}
	\end{axis}
	\end{tikzpicture}
	\caption{Computational results: ECDF for instance $R.50$}
	\label{R50}
\end{figure}			

\begin{figure}
	\begin{tikzpicture}
	\begin{axis}[height= 8cm,width=12cm,xlabel=$\pi^*$, ylabel=fraction of runs with profit within $\pi^*$,ymin=0.0,ymax=1.0,xmin=930000,xmax=1900000,
	legend style={at={(0.5,-0.15)}, anchor=north,legend columns=-1}
	]
	\addplot+ [
	%smooth,
	] coordinates {
		(1693932.675, 0.000)
		(1719036.117, 0.100)
		(1727029.129, 0.120)
		(1735022.140, 0.380)
		(1743015.151, 0.540)
		(1751008.163, 0.820)
		(1759001.174, 1.000)
	};

	\addplot+ [
	%smooth,
	] coordinates {
		(1276577.581, 0.000)
		(1319280.052, 0.080)
		(1349087.799, 0.240)
		(1378895.547, 0.320)
		(1408703.295, 0.580)
		(1438511.042, 0.780)
		(1468318.790, 1.000)
	};

	\addplot+ [
	%smooth,
	] coordinates {
		(1294008.443, 0.000)
		(1338249.393, 0.128)
		(1369419.552, 0.383)
		(1400589.710, 0.638)
		(1431759.868, 0.872)
		(1462930.027, 0.979)
		(1494100.185, 1.000)
		
	};

	\addplot+ [
	%smooth,
	] coordinates {
		(935895.853, 0.000)
		(1099107.978, 0.100)
		(1252866.609, 0.200)
		(1406625.239, 0.240)
		(1560383.870, 0.280)
		(1714142.501, 0.440)
		(1867901.132, 1.000)
	};

	\addlegendentry{BLP}
	\addlegendentry{BLPR}
	\addlegendentry{FCFS}
	\addlegendentry{PKP}
	\end{axis}
	\end{tikzpicture}
	\caption{Computational results: ECDF for instance $RC.50$}
	\label{RC50}
\end{figure}
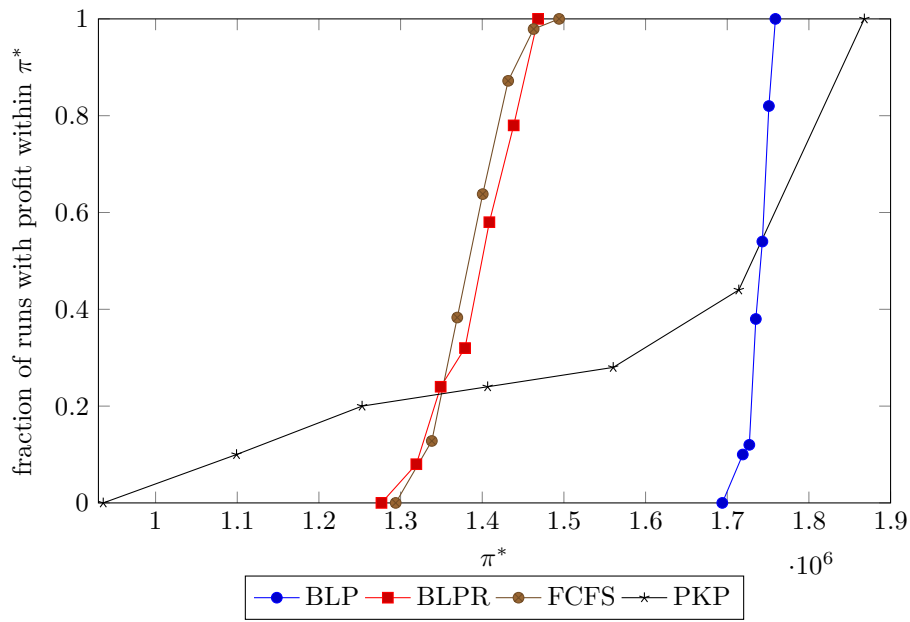
	
\newpage	
	\section{Conclusions}
	\label{section6}
	We have considered the dynamic capacity-control problem of a logistic operator that has to decide whether to accept collection requests issued by different customers on a given distribution network. The decision is aimed at balancing the current profit obtained by accepting the request, with the future expected profit deriving from preserving the resources for more valuable future requests. 
%	this being a standard decision-making problem in the  revenue management capacity-control setting. 
\color{black} We have focused on RM capacity-control techniques for distribution logistics problems where operational costs can be optimized. We have proposed a dynamic programming formulation of the problem, and have introduced a deterministic approximation along with the related approximate booking limit control policy.    In particular, for the underlying capacitated vehicle-routing problem (CVRP) considered here, we come up with a formulation that defines a novel subclass of CVRP whose aim is the profit maximization. \color{black}

%	Our main contribution is the extension to distribution logistics applications of the revenue management techniques that also account for non-negligible operational costs, the latter being evaluated by solving an optimization problem.  Furthermore, we have developed specific versions of both the deterministic approximation and the control policies, based on a distribution application formulated as a capacitated vehicle routing problem, that have been  computationally tested on a set of academic instances. 
	%\section*{References}

\end{document}